\documentclass[11pt,twocolumn]{article}

\usepackage{amsfonts}
\usepackage{graphicx}
\usepackage{flushend}

\newcommand{\vs}{\vspace{0.5ex}}

\title{A Brief Tour of Logic and Optimization}

\author{John Hooker\\Carnegie Mellon University\\ jh38@andrew.cmu.edu}

\date{November 2019\\ \vspace{5ex}}

\begin{document}
	
\maketitle

\begin{abstract}
	This paper is an informal survey of some of the deep connections between logic and optimization.  It covers George Boole's probability logic, decision diagrams, logic and cutting planes, first order predicate logic, default and nonmonotonic logics, logic and duality, and finite-domain constraint programming.  There is particular emphasis on practical optimization methods that stem from these connections, including decision-diagram based methods, logic-based Benders decomposition, and integration of CP and optimization technologies.  The paper is a slight revision of an invited article for the INFORMS Optimization Society Newsletter in observance of the 2018 Khachian Award.
\end{abstract}

\section{Introduction}

This paper offers my personal perspective on some of the deep connections between logic and optimization, based on my  research experience over the last three decades.  It gives particular attention to practical \mbox{optimization} techniques that have grown out of these connections: \mbox{decision diagrams}, logic-based Benders decomposition, and finite-domain constraint programming.  I try to keep the presentation informal but cite references that provide the details.

I am fortunate to have worked with many collaborators over the years, too many to mention here, but most of whom are cited in the references.  I would like to single out one person, however, who (though not a collaborator) was particularly inspiring to me in the early days.  That is Robert Jeroslow, a highly original thinker who transitioned to optimization \mbox{after} having studied formal logic under Anil Nerode.  This was not unlike my own background, as I wrote a PhD thesis in quantified modal logic before moving into optimization.  Aside from his ground-breaking theorem on MILP representability \cite{Jer87}, which I still much admire, he introduced an idea that I found especially insightful.  He showed that when an unsatisfiable instance of a particular type of logic problem is expressed as a linear program, the dual solution encodes a proof of unsatisfiability using a well-known inference method of formal logic.  I later realized that the dual of {\em any} optimization problem can be defined as a problem in logical inference.  This leads, among other things, to the rethinking of Benders decomposition I discuss here.  

In 1986 I attended Professor Jeroslow's series of lectures at Rutgers University \cite{Jer89}, which summed up his logic-based approach and kept me on the edge of my seat.  Unfortunately, he met an untimely death two years later, but I never forgot his enthusiasm in these lectures for the insights logic can bring to optimization.

\section{Probability Logic}

George Boole is best known for Boolean algebra, but he regarded his work in probability logic as his most important contribution \cite{Boo1854}.  It also provides a very early and remarkable link between logic and optimization.  Inference in Boole's probability logic is nothing other than a linear programming (LP) problem \cite{Hai76}.   

As an example, suppose propositions $A$, $A\rightarrow B$, and $B\rightarrow C$ have probabilities 0.9, 0.8 and 0.4, respectively.  We wish to know with what probability we can deduce $C$, if we interpret $A\rightarrow B$ as a material conditional (not-$A$ or $B$) and similarly for $B\rightarrow C$.  The 8 possible truth assignments to $A$, $B$, and $C$ are 000, 001, 010, \ldots, and have (unknown) probabilities $p_{000}, p_{001}, p_{010}, \ldots$, respectively.  The probability of $A$ is $p_{100}+p_{101}+p_{110}+p_{111}=0.9$, because $A$ is true for these 4 assignments, and similarly for the other two  propositions.  We now have LP problems that, respectively, minimize and maximize the probability of $C$:
\[
\begin{array}{ll}
\multicolumn{2}{l}{\min/\max \; p_{001}+p_{011}+p_{101}+p_{111}} \vspace{0.5ex} \\
\mbox{s.t.}& p_{100}+p_{101}+p_{110}+p_{111}=0.9 \vspace{0.5ex} \\
& p_{000}+p_{001}+p_{010}+p_{011}+p_{110}+p_{111}=0.8 \vspace{0.5ex} \\
& p_{000}+p_{001}+p_{011}+p_{100}+p_{101}+p_{111}=0.4 \vspace{1ex}  \\
& {\displaystyle \sum_i p_i = 1} \vspace{0.5ex} \\
& p_i\geq 0, \;\mbox{all}\; i
\end{array}
\]
This yields an interval $[0.1,0.4]$ of possible probabilities for $C$.  There are exponentially many variables in the LP, but it admits a straightforward application of column generation \cite{Ho88,JauHanAra91}.  
	
Similar LP models can be given for logics of belief and evidence \cite{AndHoo96}, such as Dempster-Shafer theory.  A nonlinear model can be stated for inference in a Bayesian network, where the nodes correspond to logical propositions, and the model can be simplified by exploiting the structure of the network \cite{AndHoo94}.

\section{Decision Diagrams}

Let's jump from an early example of the logic-optimization nexus to a very recent one: decision diagrams.  While there is nothing new about the concept of a decision diagram, which again derives ultimately from George Boole, its application to \mbox{optimization} and constraint programming is new on the scene.

The evolution of decision diagrams is itself an interesting story.  The world forgot about Boole's ideas for three-quarters of a century after his death, with the exception of two or three logicians.  One of those \mbox{logicians} was the unconventional philosopher Charles Sanders Pearce.  He suggested (in 1886!) that Boolean logic could be implemented in electrical switching circuits.  This idea, too, was forgotten for decades, until a young electrical engineering student at the University of Michigan was required to take a philosophy course.  I suspect that Claude Shannon, like many students today, chafed at this kind of course requirement, but it exposed him to Pearce's work. He later entered MIT (1936), where he incorporated Pearce's idea into what may be the most famous master's thesis ever written, {\em Symbolic Analysis of Relay and Switching Circuits}.  This essay not only resuscitated Boole's legacy but provided the \mbox{basis} for modern computing.

Following Shannon's work, C.\ Y.\ Lee \cite{Lee59} proposed ``binary-decision programs'' as a means of calculating the output of switching circuits.  Binary-decision programs led to the birth of binary decision diagrams (BDDs) in 1978 when S.\ B.\ \mbox{Akers} showed how to represent the programs graphically \cite{Ake78}.  BDDs remained mostly a curiosity until Randy Bryant proved in 1986 that reduced, ordered BDDs provide a unique minimal representation of a Boolean function \cite{Bry86}.  This led to widespread application of BDDs in circuit design and product configuration.

I got interested in BDDs when H.\ R.\ Andersen invited me to give a talk at IT University of Copenhagen in 2005.  There, I met his PhD student Tarik Had\v{z}i\'{c}.  Tarik and I explored BDDs and multivalued decision diagrams (MDDs) as a basis for discrete optimization and postoptimality analysis \cite{HadHoo06d,HadHoo06,HadHoo07}.  We then teamed up with Andersen and his student \mbox{P.\ Tiedemann} to propose a concept of relaxed MDDs, which we first applied to constraint programming but later became an essential tool for optimization \cite{MDDConstraintStore}.  Research on MDDs and optimization has moved in several directions since that time, including development of a general-purpose solver for discrete optimization \cite{BerCirHoeHoo14,HodvanHoeHoo10}, applications to nonlinear programming \cite{BerCir17a}, and implementations with parallel processors \cite{BerCirSabSamSarHoe12}.  Much of this research is described in a recent book \cite{BerCirHoeHoo16a} and conference \cite{DDOPT18}.  

\begin{figure}[!b]
	\centering
	\includegraphics[clip=true,trim=150 480 100 120,scale=.7]{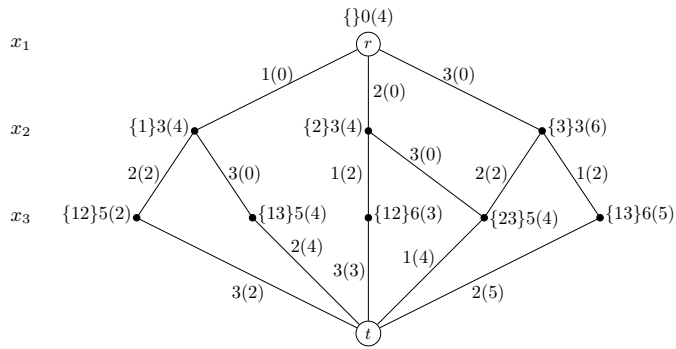}
	\vspace{-3ex}
	\caption{Decision diagram for a small job sequencing instance, with arc labels and costs shown.  States are indicated at nodes, along with minimum costs-to-go (in parentheses).}   \label{fig:DDstates}
\end{figure}

A BDD is an acyclic multigraph that represents a Boolean function.  What especially attracted me is the ability of a modified BDD or MDD  to represent the feasible set of a discrete optimization problem (by including only paths to the ``true'' terminal node). As an example, consider a job sequencing problem with time windows.  Each job $j$ begins processing no earlier than the release time $r_j$ and requires processing time $p_j$.  The objective is to minimize total tardiness, where the tardiness of job $j$ is $\max\{0,s_j+p_j-d_j\}$, and $d_j$ is the job's due date.  Figure~\ref{fig:DDstates} shows a reduced MDD for a problem instance with $(r_1,r_2,r_3)=(0,1,1)$, $(p_1,p_2,p_3)=(3,2,2)$, and $(d_1,d_2,d_3)=(5,3,5)$.  Variable $x_i$ represents the $i$th job in the sequence, arc labels indicate the value assigned to $x_i$, and arc costs appear in parentheses.  Feasible solutions correspond to paths from the root $r$ to the terminus $t$, and optimal solutions correspond to shortest paths.

The MDD is built top-down based on a dynamic programming (DP) formulation of the problem.  In the example, the DP state variables are the set of jobs already assigned and the finish time of the previous job.  The MDD is closely related to the state-transition graph for the DP, but the concept of an MDD differs in several respects.  Nodes of an MDD need not be associated with states.  MDD reduction can sometimes drastically simplify the classical state-transition graph, as in inventory management problems \cite{Hoo13}.  Finally, we will see that MDDs provide a novel basis for relaxation, primal heuristics, search, and postoptimality analysis.

\begin{figure}[!t]
	\centering
	\includegraphics[clip=true,trim=150 480 60 120,scale=.67]{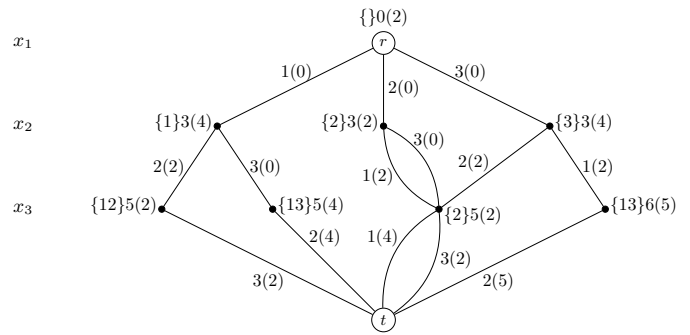}
	\vspace{-3ex}
	\caption{A relaxation of the decision diagram in Fig.~\ref{fig:DDstates}.}
	\label{fig:DDrelax}
\end{figure}

Since MDDs can grow exponentially with the problem size, much smaller relaxed MDDs are used when solving problems.  They can play a role analogous to the LP relaxation in MILP.  One way to create a relaxed MDD is to merge nodes, as illustrated in Fig.~\ref{fig:DDrelax}, which is the result of merging states $(\{1,2\},6)$ and $(\{2,3\},5)$ in layer~3 of Fig.~\ref{fig:DDstates}.  Some of the paths now correspond to infeasible solutions, including the shortest paths, which have cost 2.  This is a lower bound on the optimal cost of 4 in Fig.~\ref{fig:DDstates}.  This kind of relaxation is entirely different from the traditional state space relaxation in DP, partly because it retains the same state variables (and may require additional ones).  Sufficient conditions under which node merger results in a valid relaxation are stated in \cite{Hoo17}.  

Decision diagrams can provide all the elements of a branch-and-bound solver that does not require LP relaxation or cutting planes \cite{BerCirHoeHoo14}.
\begin{itemize}
	\item {\em Modeling framework.}  MDDs are best suited to solve problems formulated as DPs.  They provide an opportunity to exploit recursive structure as well as a novel approach to solving DPs.
	\item {\em Relaxation.}  Relaxed MDDs provide a general method of relaxing the problem without reliance on convexity, linearity, or inequality formulations.  At least in some cases, they yield tighter bounds, more rapidly, than the cutting plane technology implemented in a state-of-the-art solver \cite{BerCirHoeHoo13}.
	\item {\em Primal heuristic.} Restricted MDDs (the opposite of relaxed MDDs) provide a general-purpose primal heuristic than can be quite fast \cite{BerCirHoeYun12}.
	\item {\em Search.}  Relaxed MDDs provide a novel search method that can substantially reduce symmetry.  Rather than branch on variables, one can branch on nodes in the last exact layer of a relaxed MDD \cite{BerCirHoeHoo14}.
	\item {\em Constraint propagation.}  This lies at the heart of constraint programming (CP) solvers.  MDDs provide a more effective propagation method than the traditional domain store.  For example, they have enabled the solution of several open traveling salesman problems with time windows \cite{CirHoe13}.
	\item {\em Postoptimality analysis.}  MDDs provide an ideal tool for exploring the entire space of near-optimal solutions, since they can compactly represent this space in a form that is conveniently queried \cite{SerHoo19}.
\end{itemize}

MDDs can be combined with more traditional techniques, such as Lagrangian relaxation \cite{BerCirHoe15}.  This can allow relaxed MDDs to yield very tight bounds on combinatorial problems that are solved by heuristics, in some cases proving optimality \cite{Hoo19}.

\section{Logic and Cutting Planes}

The network of connections between logic, cutting planes, and projection is so vast and multi-faceted that I can only offer a glimpse of it here.  The pioneer in this area is H.~Paul Williams \cite{Wil77,Wil87,Wil95}.  I sought out his work when I began exploring connections between logic and integer programming three decades ago \cite{Hoo88e}.  In more recent years, I have much enjoyed collaborating with him on re-conceiving integer programming as projection \cite{WilHoo16}.  We also developed MILP formulations of distributive justice \cite{HooWil12}, an area I study in another research life.  

Let's start simple, with an illustration of the well-known resolution algorithm for propositional logic.  The first two logical clauses below jointly imply the third, which is their resolvent.
\begin{equation}
\begin{array}{r@{\;}r@{\;}r@{\;}r@{\;}r}
x_1 & \vee & x_2 & \vee & x_3 \vs \\
x_1 &      &     & \vee & \neg x_3 \vs \\
\hline
x_1 & \vee & x_2 
\end{array} \label{eq:resolution}
\end{equation}
The symbol $\vee$ is an ``or,'' $\neg$ is a ``not,'' $x_j$ and $\neg x_j$ are literals, and each variable $x_j$ can be true or false.  Two clauses can be resolved when exactly one variable $x_k$ changes sign between them.  The resolvent contains all the literals that occur in either clause except $x_k$ and $\neg x_k$.  Repeated application of resolution is a complete inference method \cite{Qui52,Qui55}.  That is, a clause set implies all and only clauses that can be obtained by a series of resolution operations (one or both parents of a resolvent may be a resolvent previously obtained).

The resolvent in (\ref{eq:resolution}) can also be viewed as a rank~1 Chv\'{a}tal-Gomory (C--G) cut by reinterpreting $x_j$ as a 0--1 variable and each clause as a linear inequality:
\begin{equation}
\begin{array}{r@{\;}r@{\;}r@{\;}r@{\;}r@{\;}lr}
x_1 & + & x_2 & + & x_3 & \geq 1 & (1/2) \vs \\
x_1 &   &     & - & x_3 & \geq 0 & (1/2) \vs \\
    &   & x_2 &   &     & \geq 0 & (1/2) \vs \\
\cline{1-6} 
\ \vspace{-2ex} \\
x_1 & + & x_2 &   &     & \geq \lceil \frac{1}{2} \rceil 
\end{array} \label{eq:resolution2}
\end{equation}
Each 1/2 in parentheses is a multiplier in the linear combination that yields the C--G cut, and bounds of the form $0\leq x_j\leq 1$ are added in when necessary to compensate for missing terms.  Chv\'{a}tal famously proved that repeated linear combination and rounding operations, where each operation may \mbox{involve} multiple inequalities, comprise a complete inference method for 0--1 inequalities \cite{Chv73}. That is, all valid \mbox{0--1} inequalities can be so obtained, a result that might be regarded as the fundamental theorem of cutting planes.  The connection with resolution is quite close. An examination of Chv\'{a}tal's elegant proof reveals that a resolution procedure lies at the heart of the argument!  

One might naturally ask if the elementary closure (set of rank 1 C--G cuts) is also somehow connected with resolution.  It is.  Input resolution is a restricted form of resolution in which at least one parent of \mbox{every} resolvent is a clause in the original set.  Clauses that can be derived by input resolution are precisely those whose inequality form belongs to the elementary closure \cite{Hoo89}.

One might also ask if a complete resolution-like procedure can be stated for general 0--1 inequalities when they are treated as logical propositions.  That is, inequalities are viewed as equivalent when they are satisfied by the same 0--1 points, and we wish to derive all valid inequalities up to logical equivalence.  The answer is again yes.  This can be accomplished by applying two types of operations repeatedly:  classical resolution and diagonal summation, where the latter derives a very particular type of rank~1 cut \cite{Hoo88,Hoo92}.  This theorem may be viewed as a logical analog of Chv\'{a}tal's cutting plane theorem.  It also provides the basis for a pseudoboolean optimization solver \cite{Bar95}.

I should remark in passing that Boolean and pseudoboolean methods comprise an entire field I am not covering here, one that has applications in both optimization and data analytics.  Any 0--1 programming problem can be viewed as optimizing a single (generally nonlinear) ``pseudoboolean'' expression, whose variables are 0--1 and whose value is a real number.  The late \mbox{Peter} \mbox{Hammer} is the seminal figure in this area \cite{HamRud68}.  In past years I made several enjoyable trips to Rutgers University to discuss ``things Boolean'' (as he put it) with him and his colleague Endre Boros, and we coauthored a couple of papers \cite{BorHamHoo94,BorHamHoo95}.  A summary of some results from the field can be found in \cite{Hoo02b}.  

Still another property of resolution is its relation to integrality.  Applying the resolution algorithm to a clause set does not ensure that its inequality form describes an integral polytope, but it reduces the question of integrality to that of underlying set covering problems \cite{Hoo96e}.

Finally, resolution can be interpreted as projection.  Returning to the example (\ref{eq:resolution}), it can be viewed as projecting out the variable $x_3$.  Its inequality form (\ref{eq:resolution2}) is an instance of Fourier-Motzkin elimination (a projection method for linear systems) plus a rounding step.  Projection is a very general idea that can form the basis of a unifying framework for logical inference, optimization, and constraint programming \cite{Hoo16,Hoo16e}.  It can suggest a more general and arguably more natural conception of integer programming (IP), in which the feasible set is conceived as a subset of the integer lattice defined by congruence relations.  Then the projection of an IP onto a subspace is again an IP, which is not the case for conventional IP.  Cutting plane theory takes a different form in which rounding is replaced by application of a generalized Chinese remainder theorem \cite{WilHoo16}.  

There is more.  In a bit of serendipity, early research on logic and IP led to the discovery of the phase transition for SAT problems \cite{HooFed90}.  An IP rounding theorem \cite{Cha84} led to a substantial extension \cite{ChaHoo91} of the class of SAT problems that can be solved by unit resolution, a linear-time inference method.  Additional research involving logic, cutting planes, and IP is described in \cite{ChaHoo99,Hoo97,Hoo00,Hooker12}.

\section{Predicate and Default Logics} 

I can only briefly mention a few of the tantalizing links between optimization and first-order predicate logic.  There is also a nice relationship with default and nonmonotonic logics.  I have fond memories of working with Vijay Chandru, Giorgio Gallo and Gabriella Rago on these topics.

First-order predicate logic, which dates back at least to Aristotle, is probably the most intensively studied form of logic.  It is characterized by such quantifiers as ``for some'' and ``for all.'' A fundamental theorem of first-order logic is Herbrand's theorem \cite{Her30}, which states roughly that an unsatisfiable formula has an unsatisfiable finite instantiation, for example a conjunction of (``ground level'') logical clauses.  

This suggests that one can check satisfiability by reducing a first-order formula to a SAT problem, perhaps \mbox{using} a method inspired by row-generation techniques in optimization.  One can in fact do this, an idea first proposed, unsurprisingly, by Robert Jeroslow \cite{Jer88a}.  Efficient ``primal'' and ``dual'' partial instantiation methods were later developed for first-order satisfiability testing \cite{ChaHoo99,Hoo93e,HooRagChaShr02}. Herbrand's theorem is itself very closely related to compactness theorems in infinite-dimensional LP and IP \cite{ChaHoo99}.

Default and other nonmonotonic logics allow an inference to be defeated when additional premises are added.  IP models for these logics not only offer a computational approach, but they actually make the underlying ideas more perspicuous than when expressed in the usual logical idiom \cite{ChaHoo99}.

\section{Logic and Duality}

Duality is one of the most beautiful and useful concepts in the field of optimization.  This seems particularly so when the connection between logic and duality is recognized.  The logical perspective also leads to a generalization of Benders decomposition that has found a variety of applications.

All optimization duals of which I am aware are, in essence, problems in logical inference.   Specifically, they can be viewed as inference duals \cite{Hoo00}.  The \mbox{inference} dual of a given problem seeks the tightest bound on the objective value that can be inferred from the constraint set, using a specified logical inference method.   

To make this more precise, consider a general optimization problem $\min \{f(x) \;|\; C(x), \; x\in D\}$, 
in which $C(x)$ represents a constraint set containing variables in $x=(x_1,\ldots, x_n)$, and $D$ is the domain of $x$ (such as tuples of nonnegative reals or integers).  
The inference dual is 
\begin{equation}
\max \Big\{ v\in \mathbb{R} \;\Big|\; C(x) \stackrel{P}{\vdash} \big(f(x) \geq v\big), \; P\in {\cal P}\Big\}
\label{eq:lbbd2}
\end{equation}
where $C(x) \stackrel{P}{\vdash} (f(x) \geq v)$ indicates that proof $P$ deduces $f(x)\geq v$ from $C(x)$ on the assumption that $x\in D$.  The domain of variable $P$ is the family $\mathcal{P}$ of proofs allowed by a specified inference method.  The solution of the dual consists of a proof $P$ that delivers the optimal bound $v$.

If the problem is an LP, the inference method traditionally consists of nonnegative linear combination and material implication.  An inequality  $g(x)\geq 0$ materially implies $h(x)\geq 0$ if any $x\in D$ that satisfies $g(x)\geq 0$ also satisfies $h(x)\geq 0$.  Thus in an LP context, a bound $cx\geq v$ can be inferred from the constraints $Ax\geq b$ when $uAx\geq ub$ materially implies $cx\geq v$ for some $u\geq 0$.  Given this inference method, (\ref{eq:lbbd2}) becomes the standard LP dual.  The solution of the inference dual is the proof encoded by optimal dual multipliers obtained from the traditional dual.  The LP dual is strong because the inference method is complete, due to the Farkas Lemma.\footnote{If the dual is expressed in the conventional fashion, it is strong (i.e, the primal and dual have the same value, possibly $\pm \infty$) only when the primal or dual is feasible.}

Other well-known duals are likewise special cases of inference duality, as shown in Table~\ref{ta:duals}.  Note the close similarity of Lagrangian and surrogate duality, whose classical formulations appear quite different.  The inference method associated with Lagrangian \mbox{duality} combines nonnegative linear combination with ``domination,'' which is only slightly different than material implication. An inequality $g(x)\geq 0$ dominates $h(x)\geq 0$ when $g(x)\leq h(x)$ for all $x\in D$, or when no $x\in D$ satisfies $g(x)\geq 0$. Since domination is stronger than material implication, it follows immediately that the Lagrangian dual never yields a tighter bound than the surrogate dual, a fact that classically is nontrivial to prove.

\begin{table}
{\small
\caption{Optimization duals as inference duals} \label{ta:duals}
\vspace{1.5ex}
\begin{tabular}{lll}
Type of dual & Inference method & Strong? \\
\hline
\ \vspace{-1.5ex} \\
LP          & Nonneg.\ linear combination & Yes$^*$ \\
            & + material implication         & \\
\ \vspace{-1.5ex} \\
Lagrangian  & Nonneg.\ linear combination & No \\
            & + domination                   & \\
\ \vspace{-1.5ex} \\
Surrogate   & Nonneg.\ linear combination & No \\
            & + material implication         & \\
\ \vspace{-1.5ex} \\
Subadditive & Cutting planes                 & Yes$^{**}$ \\
\hline
\end{tabular}
}
\vspace{-0.8ex}
\\
{\footnotesize 
$^*$Due to Farkas Lemma \vspace{-1ex} \\
$^{**}$Due to Chv\'{a}tal's theorem \cite{Chv73}
}
\end{table}


Let's now examine how inference duality offers a new perspective on Benders decomposition.  Jacques \mbox{Benders} proposed his brilliant idea for a decomposition method nearly 60 years ago \cite{Ben62}.  It has since become a standard optimization technique, with rapidly growing popularity in recent years \cite{RahCraGen17}.  \mbox{Benders} decomposition exploits the fact that a problem may radically simplify when certain variables (the master problem variables) are fixed.  Very \mbox{often}, the subproblem that remains after fixing these variables decouples into smaller problems.  Duality plays a key role in the method, because Benders cuts are based on dual multipliers obtained from the LP dual of the subproblem.   

Benders decomposition is actually a more general and more powerful method than its inventor may have realized.  The classical method is limited by the fact that the subproblem must be an LP, or a continuous nonlinear program in Geoffrion's generalization  \cite{Geoffrion72}.  This limitation can be overcome by replacing the classical dual with an inference dual, which can be defined for an arbitrary subproblem.  This results in ``logic-based'' Benders decomposition (LBBD), introduced in \cite{Hoo00,HooOtt03}.  This maneuver also reveals the core strategy of Benders-like methods: {\em the same proof} that establishes optimality of the subproblem solution (i.e., the same solution of the inference dual) {\em can prove useful bounds for other subproblems} that arise when the master problem variables are fixed to different values.  

A brief formal statement of LBBD might go as follows; more detail can be found in  \cite{Hoo06b,Hoo19e}.  LBBD is applied to a problem of the form
\begin{equation}
\min \big\{ f(x,y) \;\big|\; C(x,y), \; x\in D_x, \; y\in D_y \big\} \label{eq:prob}
\end{equation}
where $C(x,y)$ is a constraint set containing variables $x,y$.  Fixing $x$ to $\bar{x}$ 
defines the subproblem
\begin{equation}
\mathrm{SP}(\bar{x}) = \min \big\{ f(\bar{x},y) \;\big|\; C(\bar{x},y),\; y\in D_y \big\}  \label{eq:sub}
\end{equation}
where SP$(\bar{x})=\infty$ if (\ref{eq:sub}) is infeasible.  The inference dual of (\ref{eq:sub}) is
\begin{equation}
\max \Big\{ v\in \mathbb{R}\;\Big|\; C(\bar{x},y) \stackrel{P}{\vdash} \big(f(\bar{x},y) \geq v\big),\; P\in {\cal P}\Big\}  \label{eq:dual}
\end{equation}
The associated inference method is the procedure used to prove optimality of SP$(\bar{x})$, which may involve branching, cutting planes, constraint propagation, and so forth.  To ensure convergence of LBBD, we suppose the dual is a strong dual, as it is for exact optimization methods. 

Let proof $P^*$ solve (\ref{eq:dual}) by deducing the bound $f(\bar{x},y)\geq \mathrm{SP}(\bar{x})$.  The key to LBBD is that this same proof may deduce useful bounds when $x$ is fixed to values other than $\bar{x}$.  The term ``logic-based'' refers to this pivotal role of logical deduction.  A logic-based Benders cut $z\geq B_{\bar{x}}(x)$ is derived by identifying a bound $B_{\bar{x}}(x)$ that proof $P^*$ deduces for a given $x$.  Thus, in particular, $B_{\bar{x}}(\bar{x})=\mathrm{SP}(\bar{x})$.  The cut is added to the master problem, which in iteration $k$ of LBBD is
\[
z_k = \min \big\{ z \;\big|\; z \geq B_{x^i}(x),\; i=1, \ldots, k; \; x\in D_x \big\} \label{eq:master}
\]
where $x^1, \ldots, x^k$ are the solutions of the master problems in iterations $1, \ldots, k$, respectively.  The algorithm terminates when
\[
z_k=\min\{\mathrm{SP}(x^i)\;|\;i=1, \ldots, k\}
\]
or when $z_k=\infty$ (indicating an infeasible problem).  

To make this more concrete, let's apply it to a job assignment and scheduling problem \cite{Hoo06b}.  Jobs $1,\ldots,n$ must be assigned to facilities $1,\ldots, m$ and scheduled on those facilities.  Each job $j$ has processing time $p_{ij}$ on facility $i$ and release time $r_j$.  The facilities allow cumulative scheduling, meaning that jobs can run in parallel so long as the total rate of resource consumption does not exceed capacity $C_i$, where job $j$ consumes resources at rate $c_{ij}$.

The problem decomposes naturally.  If the master problem assigns jobs to processors, then the subproblem decouples into a separate scheduling problem for each facility.  We use LBBD because the subproblem is a combinatorial problem that cannot be treated with classical Benders decomposition.  We solve the master problem by MILP and the subproblem by constraint programming (CP). 

Let binary variable $x_{ij}=1$ when job $j$ is assigned to facility $i$.  If we minimize makespan $M$, the master problem (\ref{eq:prob}) is
\[
\begin{array}{l}
{\displaystyle
\min \Big\{ M \;\Big|\;M\geq M_i, \;\mbox{all $i$}; \; \sum_i x_{ij} = 1, \;\mbox{all $j$};
} \\
{\displaystyle 
\hspace{12ex}
\mbox{Benders cuts}; \; x_{ij}\in\{0,1\},\;\mbox{all $i,j$} \Big\}
}
\end{array}
\label{eq:sched1}
\]
The variable $M_i$ is the makespan on facility $i$ and will appear in the cuts.  Let $\bar{x}_{ij}$ be the solution of the master problem, and $J_i=\{i\;|\; \bar{x}_{ij}=1\}$ the set of jobs assigned to processor $i$.  If variable $s_j$ is the start time of job $j$, the subproblem for each facility $i$ can be given the CP formulation
\[
\begin{array}{l}
\min \Big\{M_i \;\Big|\; 
M_i\geq s_j+p_{ij},\; s_j\geq r_j, \; \mbox{all $j\in J_i$};  
\\
\hspace{13ex}
\mathrm{cumulative}\big(s(J_i),p_i(J_i),c_i(J_i),C_i\big)\Big\}
\end{array}
\label{eq:sched2}
\]
where $s(J_i)$ is the tuple of variables $s_j$ for $j\in J_i$, and similarly for $p_i(J_i)$ and $c_i(J_i)$.  The cumulative constraint,
a standard global constraint in CP, imposes the resource limitation.

Benders cuts can be obtained as follows.  Let $M_i^*$ be the optimal makespan obtained for facility $i$.  We wish to obtain a Benders cut $M_i\geq B_{i\bar{x}}(x)$ for each facility $i$ that bounds the makespan for any assignment $x$, where $B_{i\bar{x}}(\bar{x})=M_i^*$. An analysis of the subproblem structure \cite{Hoo06b} yields the following logic-based Benders cut for each facility $i$:
\[
M_i \geq M^*_i - \sum_{j\in J_i} p_{ij}(1-x_{ij}) - \max_{j\in J_i} \{r_j\} + \min_{j\in J_i} \{r_j\} \label{eq:makespanAnalytic2}
\]

Experience shows that LBBD can be substantially accelerated by including a relaxation of the subproblem within the master problem.  This is a different kind of relaxation than the usual, because it is not stated in terms of the variables $y$ of the problem being relaxed; rather, it is stated in terms of the master problem variables $x$.  The following is a valid relaxation of the cumulative job scheduling problem on facility $i$:
\[
M_i\geq r_j + (1/C_i)\hspace{-1.5ex}\sum_{j'\in J(r_j)} \hspace{-1.5ex} p_{ij'}c_{ij'}x_{ij'}
\]
where $J(r_j)$ is the set of jobs with release time no earlier than $r_j$.  This inequality can be added to the master problem for each facility $i$ and each distinct release time $r_j$.  This relaxation is straightforward, but a relaxation for the minimum tardiness problem can be quite interesting \cite{Hoo06b}.

An LBBD procedure based partly on the above cuts and relaxation solves the job assignment and scheduling problem several orders of magnitude more rapidly than state-of-the-art MILP and CP solvers \cite{CirCobHoo16}.  

When the master problem takes much longer to solve than the subproblem, as is frequently the case, one may benefit by using a variation of LBBD known as branch and check.  It was introduced along with standard LBBD in \cite{Hoo00} and first tested computationally in \cite{Tho01}, which coined the name ``branch and check.''  This approach solves the master problem only once and generate Benders cuts during the branching process, normally each time it discovers a solution that is feasible for the master problem.  

One drawback of LBBD is that Benders cuts must be designed specifically for each problem class, whereas the cuts are ready-made in classical Benders decomposition.  On the other hand, LBBD provides an opportunity to exploit problem structure with well-designed cuts.  

In any event, LBBD has been applied to a wide variety of problems, ranging from home health care delivery \cite{HecHooKim19} to frequency spectrum allocation \cite{Hof18,Kid19}.  The latter is a massive problem recently solved for the U.S. Federal Communications Commission, an achievement that was recognized with the 2018 Franz Edelman Award.  Some 113 articles describing LBBD applications are cited in \cite{Hoo19e}.

\section{Constraint Programming}

Constraint programming, a relatively new field with rapidly growing applications, provides an alternative approach to constraint solving and optimization.  It has a clear connection to logic, as its roots lie in logic programming.  Logic programming gave rise to constraint logic programming, which embeds constraints in the logical formalism by extending unification (a step in the application of Herbrand's theorem) to constraint solving in general \cite{JafLas87}.  Constraint logic programming then evolved into today's constraint programming (CP) technology.  More details on CP history can be found in \cite{Hoo02b}, and an exposition of the concepts (written for the optimization community) in \cite{Hooker12}.  

A salient feature of a logic program is that one can read its statements both declaratively and procedurally.  This characteristic is preserved in a CP model, where each constraint is processed individually, much like a statement in a computer program.  Unlike the atomistic inequality constraints of MILP, CP constraints can be high-level descriptions of problem elements, such as the cumulative scheduling constraint mentioned earlier.  The constraints are processed by applying constraint-specific inference algorithms that reduce (``filter'') the variable \mbox{domains}, enabling the solver to exploit substructure in the problem.  In addition, constraint propagation links the constraints by conveying reduced domains (or simplified relaxed MDDs) from one constraint to another.  

Constraint programming can be profitably integrated with MILP, thus combining logic and optimization in yet another way.  This can be accomplished by linking the two with LBBD or column generation, or by directly incorporating \mbox{filters} and constraint propagation into a branch-and-cut solver, as demonstrated by the research solver SIMPL \cite{AroHooYun04,HooOso99,YunAroHoo10} and to some extent by SCIP \cite{AchBerKocWol08}. SIMPL takes advantage of the powerful CP modeling paradigm, a move that raises some interesting modeling issues \cite{CirHooYun16}.  \mbox{Finally}, one can use a high-level modeling system such as MiniZinc, which ``flattens'' the model into MILP, CP and SAT components as required \cite{NetStuBecBraDucTac07}.  

Perhaps the most important concept in CP is consistency.  A partial assignment (i.e., an assignment of values to a subset of variables) is consistent with a constraint set if it occurs in some solution of the set.  A constraint set is itself consistent if every partial assignment that violates no constraint in the set is consistent with the set.  (A partial assignment can violate a constraint only if it assigns a value to every variable in the constraint.)  For example, the partial assignment $(x_1,x_2)=(0,0)$ violates no constraint in the 0--1 constraint set 
\begin{equation}
\begin{array}{r@{\;}r@{\;}r@{\;}r@{\;}r@{\;}r}
x_1 &+& x_2 &+& x_3 &\geq 1 \\
    & & x_2 &-& x_3 &\geq 0
\end{array} \label{eq:example0}
\end{equation}
but is inconsistent with the set.  The constraint set itself is therefore inconsistent.

Consistency is important because any consistent constraint set can be solved by a branching algorithm without backtracking.  Since full consistency is very hard to achieve, CP uses various weaker forms of consistency, such as domain consistency, \mbox{$k$-consistency}, and arc consistency.  Each reduces backtracking to a certain degree.  CP solvers rely especially on the achievement of (partial) domain consistency for individual constraints by reducing \mbox{domains} when the constraints are processed.  

Curiously, the concept of consistency apparently never developed in the mathematical programming community, \mbox{despite} its close connection with backtracking.  Nonetheless, conventional cutting planes can achieve a certain \mbox{degree} of consistency.  For example, adding the rank~1 C-G cut $x_1+x_2\geq 1$ to (\ref{eq:example0}) achieves consistency. Thus traditional cutting planes can reduce backtracking even without an LP relaxation.  This is an almost totally unexplored area.

\begin{figure}[!t]
	\centering
	\includegraphics[scale=0.55,clip,trim=350 130 280 150]{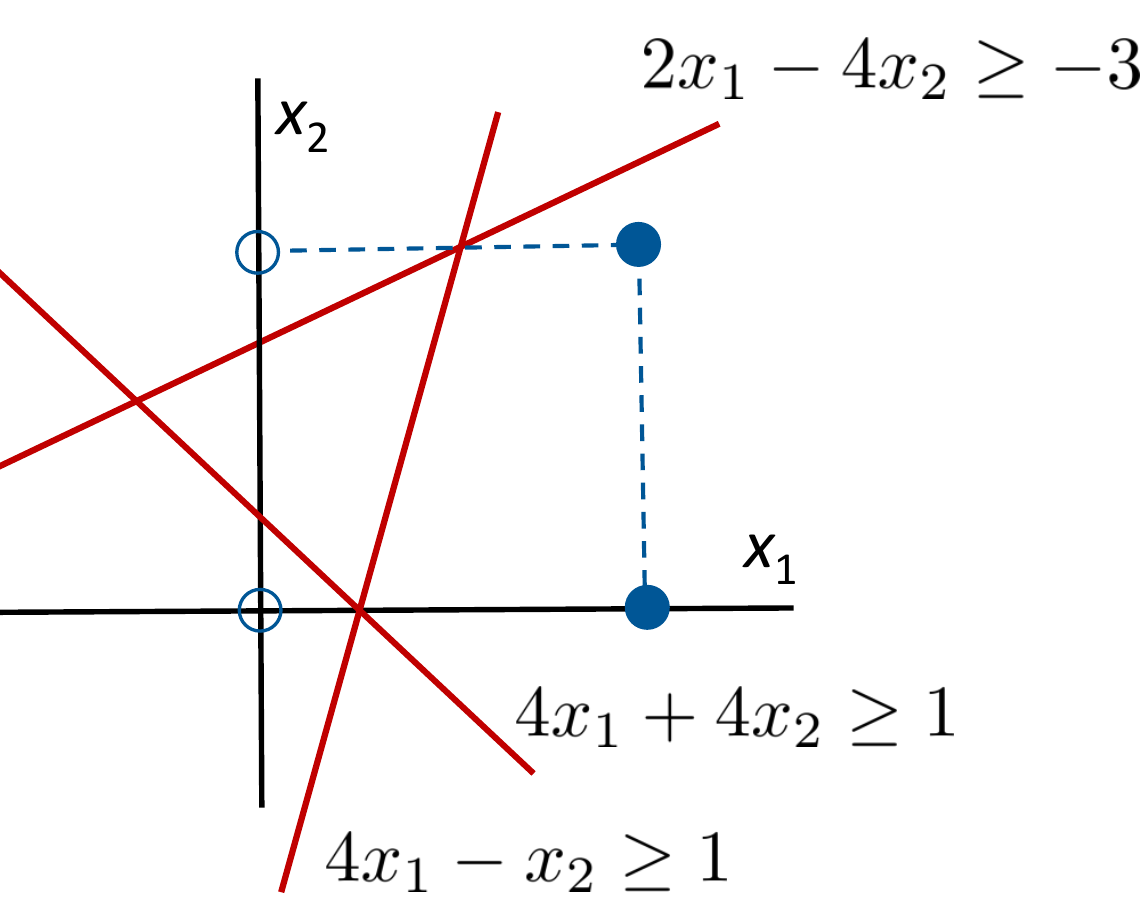} 
	\caption{Illustration of LP-consistency}\label{fig:LPconsistency3}
\end{figure}

While consistency is vary hard to achieve for an IP problem, \mbox{LP-consistency} is an attainable goal.  A \mbox{0--1} partial assignment is \mbox{LP-consistent} with constraint set $S$ if it is consistent with the LP relaxation of $S$.  Constraint set $S$ is itself \mbox{LP-consistent} if any \mbox{0--1} partial assignment that is \mbox{LP-consistent} with $S$ is consistent with $S$.  For example, the partial assignment $x_1=0$ is \mbox{LP-consistent} with the set $S$ of \mbox{0--1} constraints \mbox{$2x_1-4x_2\geq -3$} and \mbox{$4x_1-x_2\geq 1$} but is inconsistent with $S$ (Fig.~\ref{fig:LPconsistency3}).  So $S$ is not \mbox{LP-consistent}.  Adding the cut
\begin{equation}
4x_1-x_2\geq 1  \label{eq:LPcut}
\end{equation}  
achieves \mbox{LP-consistency}.  The connection with logic is preserved, because a \mbox{0--1} partial assignment is \mbox{LP-inconsistent} with $S$ if and only if it violates a logical clause whose inequality form is a rank~1 C-G cut for $S$.   Achieving \mbox{LP-consistency} avoids backtracking if one solves the LP relaxation at each node of the branching tree to check whether a variable setting is \mbox{LP-consistent}.  Thus one would not set $x_1=0$ in the example.  

Traditional cuts can achieve \mbox{LP-consistency} as well as consistency.  In the example, (\ref{eq:LPcut}) is a lift-and-project cut.  In fact, a simplified lift-and-project procedure can achieve a form of \mbox{LP-consistency} that avoids backtracking.  It is crucial to note that non-separating cuts can avoid backtracking, as the example illustrates.  If the objective is to maximize $3x_2-x_1$, the cut (\ref{eq:LPcut}) is not separating.  The  purpose of achieving \mbox{LP-consistency} is to cut off inconsistent \mbox{0--1} partial assignments, not to cut off fractional points.  The connections between logic, consistency and IP are more fully discussed in \cite{DavHoo19}.

Consistency concepts provide a new perspective on IP.  Perhaps they will form the basis for the next research thrust that combines logic and \mbox{optimization}.


\clearpage

{\small


\begin{thebibliography}{10}
	
	\bibitem{AchBerKocWol08}
	T.~Achterberg, T.~Berthold, T.~Koch, and K.~Wolter.
	\newblock A new approach to integrate {CP} and {MIP}.
	\newblock In L.~Perron and M.~A. Trick, editors, {\em CPAIOR 2008 Proceedings},
	volume 5015 of {\em Lecture Notes in Computer Science}, pages 6--20.
	Springer, 2008.
	
	\bibitem{Ake78}
	S.~B. Akers.
	\newblock Binary decision diagrams.
	\newblock {\em {IEEE} Transactions on Computers}, C-27:509--516, 1978.
	
	\bibitem{MDDConstraintStore}
	H.~R. Andersen, T.~Had\v{z}i\'{c}, J.~N. Hooker, and P.~Tiedemann.
	\newblock A constraint store based on multivalued decision diagrams.
	\newblock In {\em Principles and Practice of Constraint Programming (CP 2007)},
	volume 4741 of {\em LNCS}, pages 118--132. Springer, 2007.
	
	\bibitem{AndHoo94}
	K.~A. Andersen and J.~N. Hooker.
	\newblock Bayesian logic.
	\newblock {\em Decision Support Systems}, 11:191--210, 1994.
	
	\bibitem{AndHoo96}
	K.~A. Andersen and J.~N. Hooker.
	\newblock A linear programming framework for logics of uncertainty.
	\newblock {\em Decision Support Systems}, 16:39--53, 1996.
	
	\bibitem{AroHooYun04}
	I.~Aron, J.~N. Hooker, and T.~H. Yunes.
	\newblock {SIMPL}: A system for integrating optimization techniques.
	\newblock In J.~C. R\'{e}gin and M.~Rueher, editors, {\em CPAIOR 2004
		Proceedings}, volume 3011 of {\em Lecture Notes in Computer Science}, pages
	21--36. Springer, 2004.
	
	\bibitem{Bar95}
	P.~Barth.
	\newblock {\em Logic-based 0-1 Constraint Solving in Constraint Logic
		Programming}.
	\newblock Kluwer, Dordrecht, 1995.
	
	\bibitem{Ben62}
	J.~F. Benders.
	\newblock Partitioning procedures for solving mixed-variables programming
	problems.
	\newblock {\em Numerische Mathematik}, 4:238--252, 1962.
	
	\bibitem{BerCirSabSamSarHoe12}
	D.~Bergman, A.~Cire, A.~Sabharwal, H.~Samulowitz, V.~Sarswat, and W.-J. van
	Hoeve.
	\newblock Parallel combinatorial optimization with decision diagrams.
	\newblock In {\em CPAIOR 2012 Proceedings}, volume 8451 of {\em LNCS}, pages
	351--367. Springer, 2014.
	
	\bibitem{BerCir17a}
	D.~Bergman and A.~A Cire.
	\newblock {Discrete nonlinear optimization by state-space decompositions}.
	\newblock {\em Management Science}, 64(10):4700--4720, 2017.
	
	\bibitem{BerCirHoe15}
	D.~Bergman, A.~A. Cire, and W.-J. van Hoeve.
	\newblock Lagrangian bounds from decision diagrams.
	\newblock {\em Constraints}, 20:346--361, 2015.
	
	\bibitem{BerCirHoeHoo13}
	D.~Bergman, A.~A. Cire, W.-J. van Hoeve, and J.~N. Hooker.
	\newblock Optimization bounds from binary decision diagrams.
	\newblock {\em INFORMS Journal on Computing}, 26:253--268, 2013.
	
	\bibitem{BerCirHoeHoo16a}
	D.~Bergman, A.~A. Cire, W.-J. van Hoeve, and J.~N. Hooker.
	\newblock {\em Decision Diagrams for Optimization}.
	\newblock Springer, 2016.
	
	\bibitem{BerCirHoeHoo14}
	D.~Bergman, A.~A. Cire, W.-J. van Hoeve, and J.~N. Hooker.
	\newblock Discrete optimization with binary decision diagrams.
	\newblock {\em INFORMS Jorunal on Computing}, 28:47--66, 2016.
	
	\bibitem{BerCirHoeYun12}
	D.~Bergman, A.~A. Cire, W.-J. van Hoeve, and T.~Yunes.
	\newblock {BDD}-based heuristics for binary optimization.
	\newblock submitted.
	
	\bibitem{Boo1854}
	G.~Boole.
	\newblock {\em An Investigation of the Laws of Thought, On Which are Founded
		the Mathematical Theories of Logic and Probabilities}.
	\newblock Walton and Maberly, London, 1854.
	
	\bibitem{BorHamHoo94}
	E.~Boros, P.~Hammer, and J.~N. Hooker.
	\newblock Predicting cause-effect relationships from incomplete discrete
	observations.
	\newblock {\em SIAM Journal on Discrete Mathematics}, 7:531--543, 1994.
	
	\bibitem{BorHamHoo95}
	E.~Boros, P.~Hammer, and J.~N. Hooker.
	\newblock Boolean regression.
	\newblock {\em Annals of Operations Research}, 58:201--226, 1995.
	
	\bibitem{Bry86}
	R.~E. Bryant.
	\newblock Graph-based algorithms for boolean function manipulation.
	\newblock {\em {IEEE} Transactions on Computers}, C-35:677--691, 1986.
	
	\bibitem{Cha84}
	R.~Chandrasekaran.
	\newblock Integer programming problems for which a simple rounding type of
	algorithm works.
	\newblock In W.~R. Pulleyblank, editor, {\em Progress in Combinatorial
		Optimizaion}, pages 101--106. Academic Press Canada, 1984.
	
	\bibitem{ChaHoo91}
	V.~Chandru and J.~N. Hooker.
	\newblock Extended {Horn} clauses in propositional logic.
	\newblock {\em Journal of the ACM}, 38:203--221, 1991.
	
	\bibitem{ChaHoo99}
	V.~Chandru and J.~N. Hooker.
	\newblock {\em Optimization Methods for Logical Inference}.
	\newblock Wiley, 1999.
	
	\bibitem{Chv73}
	V.~Chv\'{a}tal.
	\newblock Edmonds polytopes and a hierarchy of combinatorial problems.
	\newblock {\em Discrete Mathematics}, 4:305--337, 1973.
	
	\bibitem{CirCobHoo16}
	A.~A. {Cir\'{e}}, E.~{\c{C}oban}, and J.~N. Hooker.
	\newblock Logic-based {Benders} decomposition for planning and scheduling: A
	computational analysis.
	\newblock {\em Knowledge Engineering Review}, 31:440--451, 2016.
	
	\bibitem{CirHooYun16}
	A.~A. {Cir\'{e}}, J.~N. Hooker, and T.~Yunes.
	\newblock Modeling with metaconstraints and semantic typing of variables.
	\newblock {\em INFORMS Journal on Computing}, 28:1--13, 2016.
	
	\bibitem{CirHoe13}
	A.~A. Cire and W.-J. van Hoeve.
	\newblock Multivalued decision diagrams for sequencing problems.
	\newblock {\em Operations Research}, 61:1411--1428, 2013.
	
	\bibitem{DavHoo19}
	D.~Davarnia and J.~N. Hooker.
	\newblock Consistency for {0--1} programming.
	\newblock In {L.-M.} Rousseau and K.~Stergiou, editors, {\em CPAIOR 2019
		Proceedings}, volume 8451 of {\em LNCS}, pages 351--367. Springer, 2014.
	
	\bibitem{DDOPT18}
	{DDOPT 2018}.
	\newblock Symposium on decision diagrams for optimization.
	\newblock Carnegie Mellon University,
	{https://sites.google.com/view/ddopt-2018}, 19--20 October 2018.
	
	\bibitem{Geoffrion72}
	A.~M. Geoffrion.
	\newblock Generalized {Benders} decomposition.
	\newblock {\em Journal of Optimization Theory and Applications}, 10:237--260,
	1972.
	
	\bibitem{HadHoo06d}
	T.~{Had\v{z}i\'{c}} and J.~N. Hooker.
	\newblock Discrete global optimization with binary decision diagrams.
	\newblock In {\em \mbox{GICOLAG} 2006}, Vienna, Austria, December 2006.
	
	\bibitem{HadHoo06}
	T.~{Had\v{z}i\'{c}} and J.~N. Hooker.
	\newblock Postoptimality analysis for integer programming using binary decision
	diagrams.
	\newblock Technical report, Carnegie Mellon University, 2006.
	
	\bibitem{HadHoo07}
	T.~{Had\v{z}i\'{c}} and J.~N. Hooker.
	\newblock Cost-bounded binary decision diagrams for \mbox{0--1} programming.
	\newblock In E.~Loute and L.~Wolsey, editors, {\em CPAIOR Proceedings}, volume
	4510 of {\em LNCS}, pages 84--98. Springer, 2007.
	
	\bibitem{Hai76}
	T.~Hailperin.
	\newblock {\em Boole's Logic and Probability}, volume~85 of {\em Studies in
		Logic and the Foundations of Mathematics}.
	\newblock North-Holland, Amsterdam, 1976.
	
	\bibitem{HamRud68}
	P.~L. Hammer and S.\ Rudeanu.
	\newblock {\em Boolean Methods in Operations Research and Related Areas}.
	\newblock Springer, New York, 1968.
	
	\bibitem{HecHooKim19}
	A.~Heching, J.~N. Hooker, and R.~Kimura.
	\newblock A logic-based {Benders} approach to home healthcare delivery.
	\newblock {\em Transportation Science}, pages 510--522, 2019.
	
	\bibitem{Her30}
	J.~Herbrand.
	\newblock Recherches sur la {th\'{e}orie} de la {d\'{e}monstration}.
	\newblock {\em Travaux de la {Soci\'{e}t\'{e}} des sciences et des lettres de
		Varsovie, Cl.\ III, math.-phys.}, 33:33--160, 1930.
	
	\bibitem{HodvanHoeHoo10}
	S.~Hoda, W.-J. van Hoeve, and John~N. Hooker.
	\newblock A systematic approach to {MDD}-based constraint programming.
	\newblock In {\em Principles and Practices of Constraint Programming (CP
		2010)}, Lecture Notes in Computer Science. Springer, 2010.
	
	\bibitem{Hof18}
	K.~Hoffmann.
	\newblock Using hybrid optimization algorithms for very-large graph problems
	and for small real-time problems.
	\newblock INFORMS Optimization Society Conference, plenary talk, 2018.
	
	\bibitem{Hoo88}
	J.~N. Hooker.
	\newblock Generalized resolution and cutting planes.
	\newblock {\em Annals of Operations Research}, 12:217--239, 1988.
	
	\bibitem{Ho88}
	J.~N. Hooker.
	\newblock A mathematical programming model for probabilistic logic.
	\newblock working paper 05-88-89, Graduate School of Industrial Administration,
	Carnegie Mellon University, 1988.
	
	\bibitem{Hoo88e}
	J.~N. Hooker.
	\newblock A quantitative approach to logical inference.
	\newblock {\em Decision Support Systems}, 4:45--69, 1988.
	
	\bibitem{Hoo89}
	J.~N. Hooker.
	\newblock Input proofs and rank one cutting planes.
	\newblock {\em ORSA Journal on Computing}, 1:137--145, 1989.
	
	\bibitem{Hoo92}
	J.~N. Hooker.
	\newblock Generalized resolution for 0-1 linear inequalities.
	\newblock {\em Annals of Mathematics and Artificial Intelligence}, 6:271--286,
	1992.
	
	\bibitem{Hoo93e}
	J.~N. Hooker.
	\newblock New methods for computing inferences in first order logic.
	\newblock {\em Annals of Operations Research}, pages 479--492, 1993.
	
	\bibitem{Hoo96e}
	J.~N. Hooker.
	\newblock Resolution and the integrality of satisfiability problems.
	\newblock {\em Mathematical Programming}, 74:1--10, 1996.
	
	\bibitem{Hoo97}
	J.~N. Hooker.
	\newblock Constraint satisfaction methods for generating valid cuts.
	\newblock In D.~L. Woodruff, editor, {\em Advances in Computational and
		Stochastic Optimization, Logic Programming and Heuristic Search}, pages
	1--30. Kluwer, Dordrecht, 1997.
	
	\bibitem{Hoo00}
	J.~N. Hooker.
	\newblock {\em Logic-Based Methods for Optimization: Combining Optimization and
		Constraint Satisfaction}.
	\newblock Wiley, New York, 2000.
	
	\bibitem{Hoo02b}
	J.~N. Hooker.
	\newblock Logic, optimization and constraint programming.
	\newblock {\em {INFORMS} Journal on Computing}, 14:295--321, 2002.
	
	\bibitem{Hoo06b}
	J.~N. Hooker.
	\newblock Planning and scheduling by logic-based {Benders} decomposition.
	\newblock {\em Operations Research}, 55:588--602, 2007.
	
	\bibitem{Hooker12}
	J.~N. Hooker.
	\newblock {\em Integrated Methods for Optimization, 2nd ed.}
	\newblock Springer, 2012.
	
	\bibitem{Hoo13}
	J.~N. Hooker.
	\newblock Decision diagrams and dynamic programming.
	\newblock In C.~Gomes and M.~Sellmann, editors, {\em CPAIOR 2013 Proceedings},
	pages 94--110., 2013.
	
	\bibitem{Hoo16}
	J.~N. Hooker.
	\newblock Projection, consistency, and {George Boole}.
	\newblock {\em Constraints}, 21:59--76, 2016.
	
	\bibitem{Hoo16e}
	J.~N. Hooker.
	\newblock Projection, inference and consistency.
	\newblock In {\em IJCAI 2016 Proceedings}, pages 4175--4179, 2016.
	
	\bibitem{Hoo17}
	J.~N. Hooker.
	\newblock Job sequencing bounds from decision diagrams.
	\newblock In J.~C. Beck, editor, {\em Principles and Practice of Constraint
		Programming (CP 2017)}, volume 10416 of {\em LNCS}, pages 565--578. Springer,
	2017.
	
	\bibitem{Hoo19}
	J.~N. Hooker.
	\newblock Improved job sequencing bounds from decision diagrams.
	\newblock In T.~Schiex and S.~de~Givry, editors, {\em Principles and Practice
		of Constraint Programming (CP 2019)}, pages 268--283. Springer, 2019.
	
	\bibitem{Hoo19e}
	J.~N. Hooker.
	\newblock Logic-based {Benders} decomposition for large-scale optimization.
	\newblock In J.~{Velasquez-Berm\'{u}dez}, M.~Khakifirooz, and M.~Fathi,
	editors, {\em Large Scale Optimization Applied to Supply Chain and Smart
		Manufacturing: Theory and Real-Life Applications}, pages 1--26. Springer,
	2019.
	
	\bibitem{HooFed90}
	J.~N. Hooker and C.~Fedjki.
	\newblock Branch-and-cut solution of inference problems in propositional logic.
	\newblock {\em Annals of Mathematics and Artificial Intelligence}, 1:123--139,
	1990.
	
	\bibitem{HooOso99}
	J.~N. Hooker and M.~A. Osorio.
	\newblock Mixed logical/linear programming.
	\newblock {\em Discrete Applied Mathematics}, 96--97:395--442, 1999.
	
	\bibitem{HooOtt03}
	J.~N. Hooker and G.~Ottosson.
	\newblock Logic-based {Benders} decomposition.
	\newblock {\em Mathematical Programming}, 96:33--60, 2003.
	
	\bibitem{HooRagChaShr02}
	J.~N. Hooker, G.~Rago, V.~Chandru, and A.~Shrivastava.
	\newblock Partial instantiation methods for inference in first order logic.
	\newblock {\em Journal of Automated Reasoning}, 28:371--396, 2002.
	
	\bibitem{HooWil12}
	J.~N. Hooker and H.~P. Williams.
	\newblock Combining equity and utilitarianism in a mathematical programming
	model.
	\newblock {\em Management Science}, 58, 2012.
	
	\bibitem{JafLas87}
	J.~Jaffar and J.-L. Lassez.
	\newblock Constraint logic programming.
	\newblock In {\em Proceedings of the 14th symposium on Principles of
		Programming Languagess}, pages 111--119, Munich, 1987.
	
	\bibitem{JauHanAra91}
	B.~Jaumard, P.~Hansen, and M.~P. {Arag\~{a}o}.
	\newblock Column generation methods for probabilistic logic.
	\newblock {\em {\mbox{INFORMS}} Journal on Computing}, 3:135--148, 1991.
	
	\bibitem{Jer87}
	R.~G. Jeroslow.
	\newblock Representability in mixed integer programming, {I: Characterization
		results}.
	\newblock {\em Discrete Applied Mathematics}, 17:223--243, 1987.
	
	\bibitem{Jer88a}
	R.~G. Jeroslow.
	\newblock Computation-oriented reductions of predicate to propositional logic.
	\newblock {\em Decision Support Systems}, 4:183--197, 1988.
	
	\bibitem{Jer89}
	R.~G. Jeroslow.
	\newblock {\em Mixed Integer Formulation}, volume~40 of {\em Logic-Based
		Decision Support}.
	\newblock North-Holland, Amsterdam, 1989.
	
	\bibitem{Kid19}
	J.~L. Kiddoo, E.~Kwerel, S.~Javid, M.~Dunford, G.~M. Epstein, C.~E. Meisch,
	K.~L. Hoffman, B.~B. Smith, A.~B. Coudert, R.~K. Sultana, J.~A. Costa,
	S.~Charbonneau, M.~Trick, I.~Segal, K.~Leyton-Brown, N.~Newman,
	A.~{Fr\'{e}chette}, D.~Menon, and P.~Salasznyk.
	\newblock Operations research enables auction to repurpose television spectrum
	for next-generation wireless technologies.
	\newblock {\em INFORMS Journal on Applied Analytics}, 49:7--22, 2019.
	
	\bibitem{Lee59}
	C.~Y. Lee.
	\newblock Representation of switching circuits by binary-decision programs.
	\newblock {\em Bell Systems Technical Journal}, 38:985--999, 1959.
	
	\bibitem{NetStuBecBraDucTac07}
	N.~Nethercote, P.~J. Stuckey, R.~Becket, S.~Brand, G.~J. Duck, and G.~Tack.
	\newblock Minizinc: {Towards} a standard {CP} modelling language.
	\newblock In C.~Bessiere, editor, {\em Principles and Practice of Constraint
		Programming (CP 2007)}, volume 4741 of {\em Lecture Notes in Computer
		Science}, pages 529--543. Springer, 2007.
	
	\bibitem{Qui52}
	W.~V. Quine.
	\newblock The problem of simplifying truth functions.
	\newblock {\em American Mathematical Monthly}, 59:521--531, 1952.
	
	\bibitem{Qui55}
	W.~V. Quine.
	\newblock A way to simplify truth functions.
	\newblock {\em American Mathematical Monthly}, 62:627--631, 1955.
	
	\bibitem{RahCraGen17}
	R.~Rahmaniani, T.~G. Crainic, M.~Gendreau, and W.~Rei.
	\newblock The {Benders} decomposition algorithm: {A} literature review.
	\newblock {\em European Journal of Operational Research}, 259:801--817, 2017.
	
	\bibitem{SerHoo19}
	T.~Serra and J.~N. Hooker.
	\newblock Compact representation of near-optimal integer programming solutions.
	\newblock {\em Mathematical Programming}, published online 2019.
	
	\bibitem{Tho01}
	E.~Thorsteinsson.
	\newblock Branch and check: A hybrid framework integrating mixed integer
	programming and constraint logic programming.
	\newblock In T.~Walsh, editor, {\em Principles and Practice of Constraint
		Programming (CP 2001)}, volume 2239 of {\em Lecture Notes in Computer
		Science}, pages 16--30. Springer, 2001.
	
	\bibitem{Wil77}
	H.~P. Williams.
	\newblock Logical problems and integer programming.
	\newblock {\em Bulletin of the Institute of Mathematics and its Implications},
	13:18--20, 1977.
	
	\bibitem{Wil87}
	H.~P. Williams.
	\newblock Linear and integer programming applied to the propositional calculus.
	\newblock {\em International Journal of Systems Research and Information
		Science}, 2:81--100, 1987.
	
	\bibitem{Wil95}
	H.~P. Williams.
	\newblock Logic applied to integer programming and integer programming applied
	to logic.
	\newblock {\em European Journal of Operations Research}, 81:605--616, 1995.
	
	\bibitem{WilHoo16}
	H.~P. Williams and J.~N. Hooker.
	\newblock Integer programming as projection.
	\newblock {\em Discrete Optimization}, 22, 2016.
	
	\bibitem{YunAroHoo10}
	T.~H. Yunes, I.~Aron, and J.~N. Hooker.
	\newblock An integrated solver for optimization problems.
	\newblock {\em Operations Research}, 58:342--356, 2010.
	
\end{thebibliography}

}
 
\end{document}